# On Theodorus' lesson in the *Theaetetus* 147d-e

## by S. Negrepontis and G. Tassopoulos

*In the celebrated* Theaetetus *147d3-e1 passage Theodorus is giving a lesson to Theaetetus and his companion, proving certain quadratic incommensurabilities. The dominant view on this passage, expressed independently by H. Cherniss and M. Burnyeat, is that Plato has no interest to, and does not, inform the reader about Theodorus' method of incommensurability proofs employed in his lesson, and that, therefore, there is no way to decide from the Platonic text, our sole information on the lesson, whether the method is anthyphairetic as Zeuthen, Becker, van der Waerden suggest, or not anthyphairetic as Hardy-Wright, Knorr maintain.*

*Knorr even claims that it would be impossible for Theodorus to possess a rigorous proof of Proposition X.2 of the* Elements, *necessary for an anthyphairetic method, because of the reliance of the proof of this Proposition, in the* Elements, *on Eudoxus principle. But it is not difficult to see that Knorr's claim is undermined by the fact that an elementary proof of Proposition X.2, in no way involving Eudoxus principle, is readily obtained either from the proof of Proposition X.3 of the* Elements *itself (which is contrapositive and thus logically equivalent to X.2), or from an ancient argument going back to the originator of the theory of proportion reported by Aristotle in the* Topics *158a-159b passage.*

*The main part of the paper is devoted to a discussion of the evidence for a distributive rendering of the plural 'hai dunameis' in the crucial statement, 'because the powers ('hai dunameis') were shown to be infinite in multitude' in 147d7-8 (contrary to the Burnyeat collective, set-theoretic rendering), necessarily resulting in an anthyphairetic method for Theodorus' proofs of incommensurability (contrary to the Cherniss-Burnyeat neutrality thesis). The crucial sentence is transitional from (a) Theodorus' contribution to (b) Theaetetus and his companion's contribution, and the distributive plural follows, in two independent ways, from the examination of each of these connections. Thus, for (a) the distributive use of 'hai dunameis' is suggested (i) in the* Theaetetus *147d3-e1 passage itself on Theodorus' lesson, by the obviously distributive use of the same word just before (147d4-5) but also after (148b1-2) the crucial statement, (ii) in the* Scholia in Platonem *by the relevant comments to the* Theaetetus *162e passage, and, with greater detail and clarity, (iii) in the* Anonymi Commentarius In Platonis Theaetetum, *by specifying the type of infinity for powers, that arose in Theodorus' proofs, to be an infinity obtained, distributively, by division within each power separately (and not by collecting all powers in a family); and, for (b) the distributive use follows from the comments of the same* Anonymi Commentarius, *according to which the idea of Theaetetus and his companion consisted in circumscribing the infinity for powers that arose in Theodorus' proofs, by passing from that infinite over to numbers and the finite, a passage that, in concrete terms, is from incommensurability over to commensurability obtained by squaring, and thus again occurring, distributively, for every power separately.*

*1. Reconstructions and neutrality for Theodorus' proofs of incommensurability*

According to the *Theaetetus* 147d3-e1 passage, Theodorus, an elder distinguished mathematician, is showing to the two precocious boys Theaetetus and young Socrates the incommensurabiity of the line segments a and r, where r is of lengrh one foot and $a^2=Cr^2$ for C=3, 5,..., 17, C non-square, stopping at 17. The concept of incommensurability ('assumetria') in the *Theaetetus* is the same as that given in Book X of the *Elements*:

***Definitions*** (X.1, 2 of the *Elements*). Let a, b be two magnitudes with a>b; we say that a, b are *commensurable* if there are a magnitude c and numbers n, m, such that a=mc, b=nc, otherwise a, b are *incommensurable*.

The reader will not fail to notice that Theodorus starts from C=3, omitting the case C=2, well known by that time, evidently due to the Pythagoreans, as Burnyeat notes[1]. The omission may suggest that Theodorus' proofs are originally due to him and not previously known. The basic question concerns the method of proof employed by Theodorus. As Hardy and Wright write: 'The question how Theodorus proved his theorems has exercised the ingenuity of every historian.'[2]. Since the only source of Theodorus' mathematical achievements is essentially the short *Theaetetus* passage under consideration, our only hope to obtain some knowledge on the method employed can only be contained in the text of the passage. It is generally agreed that any proof of incommensurabilities of the quadratic irrationals falls under one of two basic alternatives.

*1.1. The problem with Knorr's reconstruction*

The first of these alternatives is a generalisation of 'the classic indirect proof by which it was shown that the diagonal of a square must be incommensurable with its side… otherwise the same number will be both even and odd'[3], a method already known to Aristotle[4]; we will refer to this method, or any of its variations, as the *traditional* one. However, as explicitly stated in the *Theaetetus* passage, Zeuthen correctly has stipulated that the method of Theodorus 'must have been of such a kind that the application of it to each surd required to be set out separately in consequence of the variations in the numbers entering into the proofs'. Ingenuous modifications of the traditional method, so as to conform to the case by case condition set out by Zeuthen[5], have been proposed by von Fritz[6], Hardy and Wright[7], and Knorr[8].

A feature of Knorr's reconstruction[9] is that it rests on the assumption that Theodorus had some difficulty to prove the incommensurability for C=17 and thus his last proof was for the case C=15. Knorr then provides[10] a reconstruction, based on Pythagorean triads and the arithmetic of relative primes, which somehow manages to present a problem at C=17, by dealing with the cases C=4n+3 (in particular C=3,7,11,15),
C=4(2k+1)+1=8k+5 (in particular C=5,13),
C=4n+2 (in particular C=(2), 6,10,14),
C=4n (in particular C=8,12),
while devising a difficulty for the first case of non-square C=4(2k)+1=8k+1, namely for C=17.

---

[1] [Burn1], p. 502-503
[2] [HW], p.42
[3] [Burn1], p. 504
[4] e.g. in *Analutica Hustera* 41a26-31
[5] [Ze]
[6] [vF]
[7] [HR] (pp.41-43)
[8] [Kn]
[9] argued in detail in [Kn], p. 81-83
[10] in p. 181-193

Knorr's reconstruction relies crucially on the interpretation of the word 'enescheto' (occurring in 147d6) with the meaning 'entangled in something', rather than simply as 'stopped' or 'came to a standstill'; in support of such an interpretation Knorr cites Hackforth' s Notes[11]. But whether Theodorus simply 'stopped' or became 'entangled' with some difficulty, he could have done so equally at C=17 or C=19, so this particular argument does not seem to be favoring one or the other case. In fact the *Anonymous Scholiast to the Theaetetus* does not agree with Knorr, as he uses three times, instead of 'enescheto', the words 'heste' or 'stenai' (namely a simple stop)[12]. Arguments against Knorr's thesis, that Theodorus became entangled at the case C=17 not completing it, have been given by Burnyeat[13]. A further simple argument toward the same conclusion follows from the use, in the *Theaetetus* 147d6, of the word 'mechri', whose unequivocal meaning is 'as far as', 'up to and including' (in this case, the case C=17). Restricting ourselves to Plato, a sample of the occurrences of 'mechri' indicates that the term is used in the inclusive sense only[14].

*1.2. The anthyphairetic reconstructions*

The second alternative is an anthyphairetic method of proof, first proposed by Zeuthen[15], and followed, among others, by Becker[16], Mugler[17], van der Waerden,[18], Fowler[19], Kahane[20] [Ka]. In order to describe the method, we recall the definition of anthyphairesis (implied in Books VII and X of the *Elements*):

---

[11] [Ha]

[12] AST 34,35; 35,16; 35,21.

[13] [Burn1], p. 512-513

[14] The other four occurrences of 'mechri' in the *Theaetetus* confirm our claim: 'I saw Theaetetus as far as Erineus' (143b1-2);'come with me at any rate until we see' (169a2)'; 'Theod….But not further…beyond that I shall not be able to offer myself. So. It will do if you will go with me so far.' (169c6-8); 'if he were to stick up his head from below as far as the neck just here where we are' (171c1-2). (translated by M.J. Levett, revised by M. Burnyeat in [Burn2]).
Additional occurrences in Plato confirm our claim: *Hippias Major* [dub.] 281b5-c8 ('mechri Anaxagorou', clearly including Anaxagoras); *Menexenus* 239c7-240a2 ('o men protos Kutos…mechri Aiguptou erxen… tritos de Dareios… mechri Skuthon ten archen horisato', clearly incuding Egypt and the Scythians, respectively)), 242d1-4; *Phaedo*109a9-b2, 112e1-3 ('mechri tou mesou kathienai, pera d' ou'); *Cratylus* 412c8-d1 ('mechri men tou homologeisthai…epeita de amphisbetesthai'), 412e2-413b3; *Sumposion* 223d3-5 ('mechri…epeita…'); *Politeia* 423b9-c1 ('mechri toutou…pera de me'), 460e4-461a2 ('mechri tettarakontaetidos,…mechri pentekaipentekontaetous', clearly in inclusive sense in both cases), 511b3-c2 ('mechri tou anhupothetou', clearly arriving at it); *Parmenides* 153c3-5; *Sophistes* 222a2-8 ('mechri men…entautha ho sophistes kai o aspalieutes hama…poreuesthon …ektrepesthon de..apo tes zootherikes'), 267a4-6; *Phaidrus* 248c3-5; *Timaeuss* 67e6-68b1 ('mechri ton ommaton', clearly including the eyes); *Nomoi* 744d8-745a3 ('mechri tetraplasiou…pleiona de..'), 756a4-b6 ('mechri duoin…to de triton…'); 771c1-4 ('mechri ton dodeka apo mias…plen endekados', expressly including twelve), 811c6-10 ('mechri deuro', clearly inclusive), 814d7-e4 ('tes men peri palaistran dunameos to mechri deur'…peri de tes alles kineseos', in a clearly inclusive manner), 855a7-b2 ('mechri tosoutou…, to de pleon me').
Especially clear and conclusive are the several occurrences of the type 'till x, but not beyond x'. We may conclude that if Knorr is right in his interpretation of 'mechri' in the *Theaetetus* 147d6, then this meaning of 'mechri' will be unique in all of Platonic writings.

[15] [Ze]

[16] [Be]

[17] [Mu], p xx-xxv, 174-249

[18] [vdW]

[19] [F]

[20] [Ka]

*Definition*[21]. Let a, b be two magnitudes (line segments, areas, volumes), with a>b; the *anthyphairesis* of a to b is the following sequence of mutual divisions:

$a = I_0 b + e_1$, with $b > e_1$,
$b = I_1 e_1 + e_2$, with $e_1 > e_2$,
…
$e_{n-1} = I_n e_n + e_{n+1}$, with $e_n > e_{n+1}$,
$e_n = I_{n+1} e_{n+1} + e_{n+2}$, with $e_{n+1} > e_{n+1}$,
…

and we set $Anth(a, b) = [I_0, I_1, \ldots, I_n, I_{n+1}, \ldots]$ for the sequence of successive quotients of the anthyphairesis of a to b.

In case a, b are natural numbers, this process, the 'Euclidean algorithm', is always *finite*, i.e. it ends after a finite number of successive remainders, and the last (positive) remainder divides (measures) the previous remainder, and is in fact the greatest common divisor (common measure) of a, b (Propositions VII.1,2 of the *Elements*). An anthyphairesis is *infinite* if the sequence of successive *remainders* is infinite. The fundamental anthyphairetic connection between incommensurability and infinity is contained in the following

*Proposition (X.2 of the Elements).*[22] If a, b are two magnitudes, with a>b, and the anthyphairesis of a to b is infinite, then a, b are incommensurable.

---

[21] The modern version of anthyphairesis is the theory of continued fractions of real numbers (>1). An exposition of the basics of the modern theory can be found in Hardy and Wright [HR], p. 129-203, and Fowler [F], p. 303-355.

[22] The precise statement of *Proposition X.2* of the *Elements* is the following:
*'If, when the less of two unequal magnitudes is continually subtracted in turn ('anthyphairated') from the greater that which is left ('to kataleipomenon') never measures the one before it, then the two magnitudes are incommensurable.'*
It is important to understand the exact meaning of remainder ('kataleipomenon') in the statement. Let the two magnitudes in question be a, b, with a>b. At the first stage of subtraction, we subtract the smaller b from the greater a; is the 'kataleipomenon' a-b? No, that which is left will be 'kataleipomenon' only if it 'leipetai heautou [[of b]] elasson', namely only if the remainder is smaller than b. Thus if a=kb+c, and c<b, then c=a-kb is the 'kataleipomenon', but the intermediate differences a-mb, for m=1,…,k-1, do not count as anthyphairetic remainders.
Thus the assumption of Proposition X.2:
*'If, when the less of two unequal magnitudes is continually subtracted in turn ('anthyphairated') from the greater that which is left ('to kataleipomenon') never measures the one before it'*
is precisely stated as follows:
*'Let a, b are two magnitudes, for which there are two infinite sequences, a sequence $k_1, k_2, k_3,\ldots$ of natural numbers, and a sequence $c_1 > c_2 > c_3 > \ldots$ of magnitudes, such that*
$a = k_1 b + c_1$, $c_1 < b$,
$b = k_2 c_1 + c_2$, $c_2 < c_1$,
$c_1 = k_3 c_2 + c_3$, $c_3 < c_2$,
…
*ad infinitum'.*

Thus, the notion of (anthyphairetic) remainder that is meant by Euclid is contained in the following
*Definition.* Let a>b be two (homogeneous) magnitudes. If a=kb+c, with k a natural number and b>c, we will call c *the (anthyphairetic) remainder* ('kataleipomenon') of the division of a by b.
In particular the differences a-b, a-2b,…,a-(k-1)b do *not* qualify as (anthyphairetic) remainders.
Euclid nowhere mentions the term 'infinite' in connection with anthyphairesis, but his notion of infinity, implied by the sentence 'that which is left *never* measures the one before it', is understood by The argument, set out by Zeuthen, in favour of the anthyphairetic reconstruction gained initial strength because it conformed with the separate case by case proofs that is implied by the *Theaetetus* passage, but was weakened by the von Fritz and the Hardy and Wright[23] arguments, mentioned in Section 1.1.

*1.3. Critique of Knorr's argument on the impossibility of a rigorous pre-Eudoxean proof of Proposition X.2*

At this point we must examine Knorr's argument against an anthyphairetic reconstruction of Theodorus' method of proof, based not on textual considerations but on the possible chronology of Proposition X.2.

First, Knorr correctly points out[24] that Theodorus' incommensurability proofs are described by Plato in the *Theaetetus* as strict and rigorous, and not simply making these incommensurabilities seem plausible. Knorr is basing his argument on the rigorousness of Theodorus proofs on the interpretation of the word 'apephene'. Knorr, invoking the *Parmenides* 128e5-130a2 passage, showed convincingly that Plato, at least at the time that he wrote the *Theaetetus*, was using 'apophainein' and 'apodeiknunai' (proving) interchangeably. But in fact Plato himself, as Knorr notes, makes sure in the *Theaetetus* 162e that Theodorus method involved nothing less than rigorous mathematical proof. In addition, Plato's account of Theodorus proofs in the *Theaetetus* was written several decades after their actual occurrence, enough time for discovering any weaknesses in his method. We should, thus, have no doubt that Knorr is correct on this point and Theodorus proofs are indeed described by Plato as, and in fact are, rigorous.

Secondly, Knorr again correctly points out[25] that Euclid's proof of Proposition X.2 of the *Elements* makes use of the Eudoxean definition 4 of Book V, via its use of Proposition X.1, and that Eudoxus' condition was certainly not yet enunciated on Theodorus' time, unless perhaps in some intuitive manner.

Knorr then argues that, combining these two points, it must be concluded that Theodorus, in his admittedly rigorous incommensurability proofs, could not have used a non-rigorous version of Proposition X.2, and therefore his method could not be anthyphairetic. But, as we will now explain, Knorr's argument does not hold out on closer inspection.

---

the following
*Definition.* An anthyphairesis will be called *finite (or infinite)* if it generates a finite (resp. infinite) sequence of successive (anthyphairetic) remainders.
With these definitions *Proposition X.2* can now be restated as follows:
If the anthyphairesis of the magnitudes a, b, with a>b, is infinite, then a, b are incommensurable.
[23] The Hardy-Wright reconstruction assumes that Theodorus was in possession of Propositions VII.1 & 2, but lacked a general divisibility criterion such as Proposition VII.30 of the *Elements*: if a prime number p divides the product of two numbers mn, then either p divides m or p divides n, and it was for this reason that he was forced to proceed case by case with the proofs of incommensurability.
[24] [Kn], p. 75-78
[25] [Kn], p. 122-123, 256, 274, 301

First, as already explained, Euclid, in the *Elements*, does give, in Proposition X.2, a *non-elementary* demonstration (depending on Eudoxus principle, via Proposition X.1), of the statement

[X.2] *if a, b are magnitudes, a>b, with infinite anthyphairesis, then a, b are incommensurable,*

It is also true that Euclid, in the course of the very next Proposition X.3, gives a completely rigorous and *elementary* demonstration (not depending in any way on Eudoxus principle) of the statement (call it [X.3*]):

[X.3*] *if a, b are commensurable magnitudes, and a>b, then the anthyphairesis of a to b is finite,*

Knorr does not seem to realise that Proposition [X.2] is the exact contrapositive[26] to, and thus logically equivalent of, the elementary [X.3*].

It is thus seen that Euclid in two successive Propositions (X.2 and X.3*) presents two different demonstrations of *the same logically equivalent statement,* one non-elementary depending on Eudoxus condition, the other elementary, mimicking the arithmetical Proposition VII.2.

Furthermore, there is good evidence that an elementary proof of Proposition X.2, in effect already contained in the *Elements*, is old, and goes at least as far back as Theaetetus. To see this we consider an

*Alternative proof* of Proposition X.2.

Denote the magnitudes with a, b, and we assume that a>b, and a, b commensurable. So there are a line segment c and natural numbers k, m, such that a=mc, b=nc. By Propositions VII. 1 & 2, the two natural numbers m>n have a finite anthyphairesis, say

$m=k_1n+p_1$, $p_1<n$,
$n=k_2p_1+p_2$, $p_2<p_1$,
...
$p_{n-1}=k_{n+1}p_n+p_{n+1}$, $p_{n+1}<p_n$,
$p_n=k_{n+2}p_{n+1}$.

But then,

$mc=k_1nc+p_1c$, $p_1c<nc$,
$nc=k_2p_1c+p_2,c$  $p_2c<p_1c$,
...
$p_{n-1}c=k_{n+1}p_nc+p_{n+1},c$, $p_{n+1}c<p_nc$,
$p_nc=k_{n+2}p_{n+1}c$,

namely the anthyphairesis of mc=a to nc=b is finite, contradiction.

The Alternative proof of Proposition X.2 in no way employs Eudoxus' principle; it follows an idea, mentioned by Aristotle in *Topics* 158a-159b, and goes at least as back as the unnamed originator of the theory of proportion of magnitudes based on the equality of anthyphairesis. This theory employs as definition of proportion a/b=c/d the condition Anth(a,b)=Anth(c,d), and contains as an obvious, according to Aristotle, Proposition the statement:
if a,b,c are lines then a/b=ac/bc
(obvious, since it is obvious that
Anth(a,b)=Anth(ac,bc (*)).
But (*) is proved in precisely the same way that in the Alternative proof of Proposition X.2 it was proved that
Anth(a,b)= Anth(mc,nc)=Anth(m,n) (**);
in fact the Alternative proof employs a numerical version of the *Topics* Proposition.
But the originator of the anthyphairetic theory of proportion is generally acknowledged to be Theaetetus[27], hence the idea of the Alternative proof of Proposition X.2 was certainly known, and considered obvious, to Theaetetus. The obstacle concerning the proof of X.2 envisaged by Knorr simply does not exist. There is no reason, from anything we know, to exclude that Theodorus, and even the Pythagoreans, were employing Proposition X.2 as a criterion of incommensurability.

---

[26] The *contrapositive* of an implication 'p implies q' is the implication 'not q implies not p'.

Even though Euclid is admittedly confused[28] on the relevance of Eudoxus principle on the relation of infinite anthyphairesis with incommensurability, we should not conclude from what has been said that there is no such relevance. This relevance, which, as we saw, certainly does not arise with Proposition X.2, as Euclid seems to believe, does arise with the Proposition Converse to it, stating that:
[Converse X.2] *if a, b are incommensurable magnitudes, a>b, then the anthyphairesis of a to b is infinite*[29].
or with the, logically equivalent, contrapositive Proposition stating that:
[Converse X.3*] *if a, b are magnitudes, with a>b, with finite anthyphairesis, then a,b are commensurable,*
Surprisingly Euclid does not establish anywhere in his *Elements* the [Converse X.2] Proposition (or its equivalent contrapositive [Converse X.3*]), even though it is **the only one** that is in essential need of Eudoxus principle; instead he has needlessly employed Eudoxus principle for Proposition X.2 itself.

We conclude: the fact that Euclid is in a curious confusion
(a) presenting an elementary, but perfectly rigorous proof of a proposition (X.3*) and a non-elementary (depending on Eudoxus principle) of its (logically equivalent) contrapositive proposition (X.2), and
(b) omitting any mention of the converses of these propositions (X.2 & X.3*), which do depend in an essential way on Eudoxus principle,
in no way implies that these propositions (X.2 & X.3*) were not available to earlier geometers (Theaetetus, Theodorus, and even the Pythagoreans), and does not provide an argument in support of Knorr's thesis against the possibility of a pre-Eudoxean rigorous proof of Proposition X.2.

---

[27] [Kn], p. 298-313
[28] It is not likely that Euclid's confusion has to do with contrapositive and converse statements, given his masterly exposition of, say, the fifth postulate and its contrapositive in Proposition I.29, and of the Propositions I.27 and 28, inverse of the fifth postulate, or of Propositions X.8,7, contrapositive to Propositions 5,6, respectively).
[29] Note e.g. that if for a>b Eudoxus condition fails, then a,b are incommensurable but their anthyphairesis is certainly not infinite, according to the definition of infinite anthyphairesis explained in footnote 22, since not even a single anthyphairetic remainder is produced.

So there is really no force in Knorr's argument against a rigorous pre-Eudoxean demonstration of incommensurability based on Proposition X.2.

*1.4. The Cherniss-Burnyeat neutrality thesis and its problem with the crucial statement*

(a) *The Cherniss-Burnyeat neutrality thesis*

By 1.1 and 1.3 above, neither a modified traditional reconstruction (such as the one proposed by Hardy and Wright) (by 1.1), nor an anthyphairetic one (by 1.3) can be excluded *a priori* as a model for Theodorus' method, and no received reconstruction, either ttraditional or anthyphairetic, in any of the variations, claims to find support in the *Theaetetus* passage itself. The dominant thesis, expressed by Cherniss and Burnyeat, is that the *Theaetetus* 147d3-e1 passage is *neutral*, in that it contains no information on the method employed by Theodorus. Cherniss[30] states that 'The sole evidence for Theodorus' demonstration is the passage of the Theaetetus, and that reveals only that he began with square root of 3 and selecting one surd after another up to square root of seventeen there somehow stopped. This tells us nothing of the method that he used', while Burnyeat writes that 'We are not told how [Theodorus] proved this result'[31], and that "Plato has no motive to indicate to the reader whether he has in mind any definite method of proof…there is no good reason to expect that the answer is to be squeezed out of one ambiguous sentence in Plato's dialogue'[32]

The neutrality of the Cherniss-Burnyeat thesis possesses the advantage of minimality, of not having to explain too much, and it holds well for the main body of the passage, where Theodorus' lesson is described (147d3-6).

*(b) A heuristic discussion suggesting a problem with Burneyat's rendering of the crucial statement* 147d7-8

However the immediately succeeding statement (147d7-8), which we will call 'the crucial statement',

[b] *'since the powers ('dunameis') appeared (or were shown) ('ephainonto') to be infinite in multitude ('apeiroi to plethos'), it occurred to us,...'*,

needs careful consideration. Burnyeat, who, also calls it a' key sentence'[33], describes [b] as

'recounting the thoughts suggested to himself [[Theaetetus]] and his companion by and during Theodorus' lesson',

and thus regards [b] as transitional from an account of Theodorus' lesson, described just before [b], in

---

[30] [C], p. 410-414, answering the admittedly inconclusive arguments by Mugler [Mu].
[31] [Burn1], p. 494
[32] [Burn1], p. 505
[33] [Burn1], p. 501

[a] *'Theodorus was writing something about powers for us, about the three foot and five foot power, 'apophainon' that they are not commensurable ('ou summetroi') to the one foot line'* (147d3-5)

and renders [b] as

*'since the dunameis[34] were turning out('ephainonto') to be unlimited in number, it occurred to us...'*.

Thus Burnyeat assumes that Theodorus lesson described in [a] prompted, suggested to Theaetetus and his companion the notion that 'the dunameis were turning out to be unlimited in number' in [b]. Burnyeat does not pinpoint the precise connection between [a] and [b] (a linguistic connection that will be explained in section 2), but that [b] is prompted by [a] is a reasonable assumption, already suggested by the word 'since' ('epeide', 147d7)). The connection between [a] and [b] shows that [b] must be taken into consideration in investigating the nature of Theodorus' lesson in [a], and may well throw some light on it, confirming or rejecting the Cherniss-Burnyeat neutrality thesis.

Burnyeat considers two possible renderings for [b], namely

[b'] 'since the series of whole [[non-square]] number squares (or sides of such squares) were turning out to be unlimited in number , it occurred to us...', or

[b''] 'since the squares with incommensurable sides were turning out to be unlimited in number, it occurred to us...'.

He rejects [b'], namely that Theaetetus' recounting consisted simply in 'the idea that there is *an endless series* of whole number squares (or sides of such squares)', presumably because this is a rather trivial conclusion that 'would hardly need to be prompted by a process as protracted as Theodorus' lesson', and opts rather for [b'''], namely for the idea 'that there is an indefinite, perhaps *infinite* number of squares with incommensurable sides'. Thus, according to Burnyeat[35], [b] should not be rendered as [b'] but as [b'']. From Burnyeat's point of view, the rendering of [b] as [b''] has a distinct advantage: it conveniently disposes of the rather troublesome statement [b], by turning it into an innocuous one [b''] that does not contradict, and does not pose any menace to, the basic Cherniss-Burnyeat thesis that Plato provides no information on Theodorus' method of proof.

---

[34] [Burn1], p.495-502, especially 501-502. Burnyeat writes: 'Immense heat has been generated over the terminology of the passage… The cause of the trouble is that while dunamis is applied to incommensurable lines in the definition at 148b1, it is also…used earlier to specify that Theodorus' demonstrations were about, so it has been a matter of controversy whether at that earlier stage the term stands for the *sides* of a series of squares or for the *squares* themselves, both of which were involved in the exercise.'.
Indeed, in 147d3-6 a power is a square $a^2$ commensurable with the assumed one-square foot square $r^2$ by a relation $a^2=Cr^2$ for all non square numbers C from C=3 to C=17, with its side a, as shown by Theodorus, incommensurable in length with the one foot line r; in 148a7-b2 a power is a line a such that $a^2=Cr^2$ for some non-square number C, and in consequence, as shown by Theaetetus and his companion, a line a incommensurable with the one foot line r, with its square $a^2$ commensurable with the one-square foot square $r^2$. Thus it may be the case that Theodorus might consider a power more as in 147d3-6, while Theaetetus and his companion more as in 148a7-b2.
Burnyeat argues that the controversy arose 'through failure to make any distinction between meaning and application.': 'dunamis is *applied* to incommensurable lines without meaning 'line' or 'incommensurable line'', but is *meant* to be 'square', or 'square with incommensurable side'. In any case, as pointed out by Burnyeat, the difference in the versions of the concept of power is really only linguistic, and from a geometrical point of view, inessential. It makes no difference whether we have a commensurable square with an incommensurable side, or an incommensurable side with a commensurable square.
[35] [Burn1], p. 501-502

But there is a problem; Burnyeat evidently does not realize that the rendering [b''] is open to precisely the same criticism that [b'] is: [b''] 'would hardly need to be prompted by a process as protracted as Theodorus' lesson', in fact it would be obvious before even the start of Theodorus' lesson (in no way depending on the highly non-obvious proposition that *every* power is incommensurable to the one-foot line, evidently proved by Theaetetus and young Socrates and related to Proposition X.9 of the *Elements*), since, for any number of the form $C=2n^2$, n=1, 2, 3,…(and thus for an *infinite multitude* of non-square numbers C), the incommensurability of the side a of the squares $a^2=Cr^2=2n^2r^2$ follows immediately from the Pythagorean incommensurability of the side y of the square $y^2=2r^2$, an incommensurability that Theodorus clearly considers as known at the start of his lesson[36]. Thus [b''] was known to hold even before the start of Theodorus' lesson and cannot credibly be claimed to be a thought resulting from that lesson.

We can take this logic to its natural conclusion. Theodorus' lesson would have to prompt to Theaetetus and his companion a conjecture stronger than [b'']. Indeed, when someone is demonstrating to you that the sides of squares with 3, 5, 6,…, 17 square feet are incommensurable to the one foot line, *omitting only* the square numbers 4, 9, 16 *and no others*, then the reasonable conjecture that might be suggested to you by these demonstrations would be, not simply that '…the squares with incommensurable sides were turning out to be unlimited in number', as in [b''] (something which, as we saw, was known before Theodorus' lesson), but that '*all powers corresponding to non-square numbers were turning out to have incommensurable sides*', as in

[b*] 'since *all* powers corresponding to non-square numbers were turning out to have incommensurable sides, it occurred to us…'.

Thus the infinity that enters in [b] cannot reasonably refer simply to an infinite multitude of powers, as in the [b''] rendering, since the only reasonable thought by Theaetetus would have to involve *all powers*, as in the *universal* statement [b*],
a statement evidently equivalent to

[b**] 'since **every** power corresponding to non-square numbers was turning out to have incommensurable sides, it occurred to us…'.

But statement [b**], as it stands, although it qualifies as an improvement over [b''] for a reasonable description of Theaetetus' thoughts as he follows Theodorus' lesson, it certainly cannot qualify as a rendering of [b], simply because reference to the infinity that appears in [b] is completely missing in [b*]. The role of infinity of the transitional statement [b] would have to be a different one than the one assigned by the Burnyeat [b''] rendering, but an existent one, and we are thus led to a statement of the sort [b***] 'since **every** power corresponding to non-square number was turning out to be infinite, it occurred to us…'.

---

[36] We owe this remark to Dr. Andreas Koutsogiannis

Thus starting with the Burnyeat reading [b''] of [b], and subjecting it to two reasoned modifications, we come to statement [b***], which is both a rendering of [b] and a credible account of thoughts resulting from Theodorus lesson; but [b***] is a statement with a meaning that at least at first sight is not clear, and further is no longer an 'innocuous' one for the basic Cherniss-Burnyeat thesis, but one putting it in jeopardy.

Clarification of the meaning of the statement [b***] and strong support for a [b***] reading of [b] will come, as we shall see in section 2, from the *Theaetetus* passage itself and from the ancient anonymous scholiasts to this passage.

*2. Incommensurability and infinity*

*2.1. The anthyphairetic connection between incommensurability and infinite divisibility or infinity*

Any demonstration of incommensurability, of a magnitude with respect to another magnitude, with an anthyphairetic method must necessarily rely on the use of some form of Proposition X.2, namely must demonstrate that the process of anthyphairetic division of the pair of magnitudes is *infinite*. Thus, 'infinite divisibility', or simply 'infinity', in connection with incommensurability signifies in unambiguous terms infinity of their anthypairesis, and thus incommensurability by means of Proposition X.2. No other demonstration of incommensurability, such as the traditional in any of its forms, is even remotely connected with infinite divisibility or a form of infinity. There should then be no hesitation to ascribe an anthyphairestic method any time infinite divisibility or infinity is related to incommensurability. The only conceivable kind of infinity or of infinite divisibility related to incommensurability is that of infinite anthyphairesis, and the only conceivable connection is by means of Proposition X.2.

Various ancient commentators do indeed connect incommensurability with infinite divisibility or, more generally, with a principle of infinity. Here are some such examples.

(a) *Incommensurability related to infinity divisibility in geometry*

-Proclus, in his *Commentary to Euclid* (=*in Eucliden*)[37] expresses precisely this link:

*'the irrational[38] has a place only where infinite divisibility is possible'*.

-Similarly *Scholia in Eucliden* V.22:

---

[37] 60,15-16; [Mo], p. 48
[38] By 'irrational' is meant 'incommensurable with respect to some assumed line'.

> *'There are incommensurable magnitudes, because magnitude is divisible ad infinitum, as the diameter is incommensurable to the side.'*

-In the *Anonymous Scholion in Eucliden* X.1, lines 21-47 it is stated repeatedly that the Pythagoreans were relating 'infinite divisibility' with incommensurability and finite divisibility with commensurability.

[1] According to lines 25-27, the Pythagoreans were the first to discover the method of finding the greatest commom measure of two numbers by the process of anthyphairesis; the sentence

> *'on the one hand every and any number being divided ('diairoumenon') according to any cuts leaves ('katalimpanein') a part least and indivisible'* describes a process that always ends after a finite number of steps, because the remainders are natural numbers that form a strictly decreasing sequence, so that a remainder, part ('morion'), will eventually divide precisely, without remainder, the preceding one. This is the content of Propositions VII, 1&2 of the *Elements*, according to which every anthyphairesis between two numbers is finite.

[2] According to lines 27-31, the Pythagoreans, next, were the first to realise that the corresponding process for geometrical magnitudes need not be finite, but it may result in a magnitude being infinitely divisible:

> *'but on the other hand every magnitude being divided ('diairoumenon') ad infinitum ('ep' apeiron')...'*

Thus a magnitude may well be *'ep'apeiron diairoumenon'*, equivalently 'ep' apeiron temnomenon', and result in *'apeira moria'*, each of which is also infinitely divisible. This is the content of Proposition X.2 of the *Elements*, according to which every infinite anthyphairesis between two magnitudes results in incommensurability.

[3] According to lines 42-47, the Pythagoreans were thus led to distinguish between commensurability and incommensurability for magnitudes.

It is then clear, according to *Scholion* X.1, 21-47, that the Pythagoreans closely connected incommensurability with *infinity* and infinite divisibility,

-Several ancient comments[39] correlate incommensurability with division *ad infinitum* with the argument that (i) the infinite divisibility of a geometric magnitude results in smaller and smaller magnitudes with no end, hence in the absence of a least and indivisible unit for magnitudes, as is the case with natural numbers, and (ii) the lack of a common unit as in numbers implies incommensurability. Statement (i) is true, but it is easy to see that Statement (ii) is not true. Indeed, choose and fix a magnitude, say a, and let C be the class of all magnitudes of the form ma/n for every natural numbers m and n. It is easy to verify that the class C lacks a common unit but that every pair of elements in C is a commensurable pair, and thus presents a counterexample to Statement (ii).
Thus the only valid correlation between incommensurability and infinite divisibility in geometry is anthyphairetic, namely via Proposition X.2 of the *Elements*.

---

[39]Heron, *Definitiones* 136,34 (=*Scholia in Eucliden* X.9), Proclus, *in Eucliden* 278,19-24 (= Xenocrates, *Fragment 131*, 18-22), *Scholia in Eucliden* X.28

-Proclus, in his *Commentary to Euclid* (=*in Eucliden*)[40] connects incommensurability with *an infinite regress of decreasing gnomons,* an infinity that can only refer to anthyphairetic division[41]:

*'The statement that every ratio is irrational ('arretos') belongs to arithmic only and not at all to geometry, for geometry contains irrational ratios. Likewise the definition of gnomons of squares according to the lower limit ('kata to elasson') is peculiar to arithmetic; in geometry a least ('elachiston') gnomon has no place at all.'*

(b) *Incommensurability connected to a principle of infinity*

-According to Iamblichus[42] the Pythagoreans believed that the mathematical knowledge on *incommensurability and irrationality* they had discovered should be kept secret, and he who makes it public is to meet with divine justice by drowning in *the sea*. At first glance the reader might simply think that there is nothing special with drowning in the sea, and the choice of punishment is of no importance. Indeed, what might be the link between incommensurability and the sea? The answer is given with clarity in the Pappus' *Commentary to Book X of the Elements*[43], partly preserved in the *Anonymon Scholion in Eucliden* X.1, 70-79; according to this passage the Pythagoreans believed that anyone who would have the temerity, despite oaths of secrecy to the contrary, to make public the secret mathematical dogmas on *incommensurability* and *irrationality ('alogon'),* would be punished by divine justice by drowning into the *'infinite sea* of dissimilarity', of the *Politicus* 273d6-e1 pssage[44]. This symbolic (or real) story reveals the direct link that the Pythagoreans believed existed between *incommensuraility* and *infinity.*

-The Pythagorean principle of Infinity, specifically inspired from their association with Mathematics, according to Aristotle' *Metaphysica*[45], is clearly described, by Aristotle, in *Physica*[46], Simplicius in his *Commentary on Physica*[47], and Philoponus,

---

[40] 60,9-12
[41] The connection is suggested in detail in [N2].
[42] *De Communi Mathematica Scientia* 25,27-33, *De Vita Pythagorica* 18, 88,13-89,1; 34,246,10-247,14
[43] cf. Thomson [Thoms]
[44] Proclus, in several cases, refers to precisely this passage of the *Politicus* and in every case he identifies 'ton apeiron pontoon tes anomoiotetos' with the Pythagorean infinity and with the Platonic Philebean infinity. Some of these references by Proclus are as follows:
*eis Timaion* 1,113,29-31 (where it is identified with 'ton ponton tes geneseos', the expression in the *Anonymon Scholion* ), 1,175,18-24 (where it is stated that the Platonic 'ton apeiron pontoon tes anomoiotetos' is characterized by infinite divisibility, is identified with 'ametria', with the Pythagorean principle Apeiron as described in Aristotle's *Metaphusica*, with the Apeiron of Philolaos (*fragments* 1,2,6)), *eis Parmeniden* 1009,19-1010,26 (where it is described as a cause of division), *eis Alcibiaden* 33,17-34,10 (where it is correlated with divisibility), *eis Politeian* 2,69,17.
The second related expression that appears in the *Anonymon Scholion to the Elements* X.1, 70-79, is 'tois astatois reumasin', the unstable fluctions' is also in close relation with the Infinite in several passages by Iamblichus, Porphyrius and Proclus, as follows:
Iamblichus, *De mysteriis* 1,3,19-29 and 7,2,5-11. Porphyrius, *eis Politeian*, Fragment 92. Proclus, *eis Politeian* 2,157,19-22; 2,261,22-26; 2,264,1-4; 2,275,8-14; 2,347,17-348,16, *eis Timaion* 1,91,14-23. 1,344,26-348,7. 1,385,18-386,8. 3,329,4-331,26, *eis Alcibiaden* 21,8-24,9; 87,4-9; 92,16-93,6; 143,18-144,1; 173,12-17; *eis Parmeniden* 910,6-13; 1003,2-6, *Platonic Theology* 4,85,17-87,4.
[45] 985b23-986a26 and 989b29-990a34
[46] 204a8-34, 213b22-29
[47] 455,15-459,3 (especially 455,20-456,3)

in his *Commentary on Physica*[48], as a division *ad infinitum* of magnitudes.

-The following explicit link between the philosophical Platonic *principle of Infinity* in the *Philebus* 16c5-10, 23b-25b and geometric *incommensurability* is stated in Proclus' *Commentary to Euclid*[49]:

*'If there were no [Platonic Philebean principle of] infinity, all magnitudes would be commensurable and there would be nothing 'arrheton' or 'alogon', features that are thought to distinguish geometry from arithmetic'.*

In the passages cited incommensurability is directly related either to infinite divisibility or to a principle of the infinite. In (a) it is clear that 'infinite divisibility' is directly related to infinite anthyphairesis. In (b) incommensurability is related to a philosophical principle of infinity, either Pythagorean or Platonic-Philebean. These associations suggest that these philosophic principles of infinity are related to infinite anthyphairesis and Proposition X.2 of the *Elements*. This is not the place to expand on this matter, but it should be remarked that the link of both the Pythagorean[50] and the Philebean[51] principle of infinity with incommensurability and Proposition X.2 is strong.

Our conclusion is this: even though the term 'anthyphairesis' may not, and usually is not, explicitly mentioned, there is no other conceivable link of incommensurability with infinity in an ancient mathematical or philosophic passage, but the infinity of anthyphairesis as described in Proposition X.2 of the *Elements*.

*2.2. The connection between Theodorus' incommensurability proofs and infinite divisibility in* Scholia in Platonem

It follows from our comments in 2.1 that an association of incommensurability with infinite divisibility should be construed as a strong indication of an anthyphairetic approach to incommensurability. Such an association occurs is found in the commentary *Scholia in Platonem* (abbreviated SP) on the following *Theaetetus* 162e6-7 passage, concerning Theodorus' demonstrations:

---

[48] 386,14-395,7 (especially 388,24-389,20)

[49] 6,19-22 (translation by G.R. Morrow [Mo], p. 5)

[50] These remarks point to an anthyphairetic interpretation of the Pythagorean principle of infinity. The detailed arguments in favour of such an interpretation are given in Negrepontis [N2].

[51] The link between the indefinite dyad, namely the Philebean Infinite, and the side and diameter numbers[51], the 'convergents' closely connected with infinite anthyphairesis, was realised by Taylor[51] reading some passages of the *Epinomis* (especially 990b5-991b4). There is a clearer link of the Platonic description in the *Philebus* 23b-25b between the Platonic Philebean principle of Infinity and the geometric infinite anthyphairesis. Indeed the Philebean infinite is described as an 'indefinite dyad' of two opposites, such as the great and the small or the cold and the warm, which is infinite, an infinite process, precisely because at every stage it has the power to renew the initial opposition, while in case the dyad cannot at some stage renew this opposition, then it becomes a 'finite', characterised by 'commensurability' (25e1), exactly as in Propositions X.2 & 3 of the *Elements* finite anthyphairesis is characterised by commensurability, and by a relation 'as number to number' (25a8), exactly as in Propositions X.5 &6 of the *Elements* finite anthyphairesis is characterised by a relation 'as number to number'. Thus the description of the infinite in the *Philebus* 23b-25b fits precisely to that of geometric infinite anthyphairesis in Propositions X.2-8 of the *Elements*, and the Philebean principle of infinity is revealed to be a philosophic, dialectic version of the geometrical infinite anthyphairesis. The full arguments for an anthyphairetic interpretation of the Platonic principles of infinite and finite in the *Philebus* 23b-25b are given in Negrepontis [N1].

*'Imagine how utterly worthless Theodorus or any geometer would be if he were prepared to rely on probability to do geometry.'*

The commentary in *Scholia in Platonem* runs as follows:

*'If we accept the judgement of the many as the dominant one in geometry, we would be ridiculous to claim that maginitudes are **incommensurable** to each other, and that **the finite line is divisible ad infinitum** ('diaireten eis apeiron'), and **the like** (\ta toiauta').'*

The following reading of the Commentary is suggested: Theodorus, when proving the incommensurabilities of the powers (of 3,5,…,17 feet) with respect to the one foot line, employed, not any considerations of probability (which would be worthless and make him ridiculous), but logically impeccable proofs ('apodeixin de kai anagken', 162e4-5) involving the infinite divisibility (namely the infinite anthyphairesis) of these finite lines with respect to the one foot line.

Thus the commentator suggests that Theodorus' method of proof is anthyphairetic. It may be countered that the two statements in the *Scholia*: (i) 'magnitudes are incommensurable' and (ii) 'the finite line is infinitely divisible' are independent of each other, and that it is not in the intention of the commentator to relate infinite divisibility to incommensurability. However the anonymous comment is about Theodorus' demonstrations reported in the *Theaetetus* 162e6-7, and the only Theodorus' demonstrations mentioned in the *Theaetetus* were those reported in 147d-e, showing incommensurabilities and infinity of powers, not dividing the line *ad infinitum* for some other purpose. Also, the final term (iii) 'the like' ('ta toiauta') indicates that the two statements (i) and (ii) are closely related.

We must also note that there is nothing paradoxical or against 'the judgement of the many' and common intuition in the simple divisibility of a line, in dividing, say, each part in half *ad infinitum*; however the incommensurability (and the infinite division of a line by anthyphairesis) definitely did run against the opinion of the many (cf. Plato, *Leges* 819d-820b and Aristotle, *Metaphysics* 983a12-23).

*3. The crucial statement: connection with Theodorus' incommensurability demonstrations*

The heuristic arguments in section 1.4 cast doubt on the Cherniss and Burnyeat thesis, as long as there is no satisfactory interpretation of the crucial statement [b] consistent with that thesis. In the present section we will argue, that the *Theaetetus* 147d3-e1 passage itself (in 3.1) and the ancient comments on this passage contained in *Anonymi Commentarius In Platonis Theaetetum*) (in 3.2) lead to an interpretation of the crucial statement [b], radically different from the Burnyeat interpretation, inconsistent with the Cherniss-Burnyeat thesis, and suggesting that Plato does indicate that the nature of Theodorus' method of incommensurability proofs is by infinite division, and thus anthyphairetic, according to the comments in section 2.

*3.1. The connection between the Theaetetus' and the infinite in the* Theaetetus *147d-e passage.*

The connection between incommensurability and infinity in multitude of every power,

suggested in the *Scholia in Platonem* (Section 2.2), may actually be extracted from a careful and novel reading of the *Theaetetus* 147d-e passage itself.

*(a) The bond between 'apophainon' in [a] and 'ephainonto' in [b]*

It was noted in section 1.4 that Burnyeat, although he is basing his interpretation of [b] on the connection of [b] with [a], nevertheless he does not take into account the linguistic link with which Plato manifests this connection. Indeed the two words 'apophainon' (147d4) in [a] and 'ephainonto' (147d8) in [b], from the same root, set up a bond between the two sentences [a] and [b]. It may be countered that mere chance brought these two words in close proximity, and no significant conclusion need follow from it. But the same bond appears fairly frequently in Plato, in fact a little later in the *Theaetetus* 151d-e and also elsewhere [52], thus strengthening the suggestion that Plato had indeed the intention to express by this proximity a logical bond.

Not only Burnyeat but scholars and translators in general do not activate this link. Thus H.N. Fowler's translation renders 'apophainon' as 'showing' and 'ephainonto' as 'appeared'[53], M. J. Levett's 'apophainon' as 'showing' and 'ephainonto' as 'were turning out'[54], and W. R. Knorr's 'apophainein' as 'demonstrating' and 'ephainonto' as 'recognized'[55].

But if contact with the original text is to be preserved, the two words 'apophainon' and 'ephainonto' should be rendered with words of the same root. Which words? It was argued in section 1.3, following Knorr, that 'apophainon' shound be rendered as 'was showing', since Theodorus' proofs were mathematically rigorous[56]; on account of the link described, 'ephainonto' should then be rendered as 'were shown'. Thus, according to [a], Theodorus 'was showing' that the powers are incommensurable to the one foot line, and, according to [b], for this reason the powers 'were shown' to Theaetetus and his companion to be infinite in multitude. Thus the fact that [b] holds because [a] holds, assumed loosely by Burnyeat, indeed follows in precise form from the *Theaetetus* text. But the precise link between [a] and [b] has consequences

---

[52] Cf. *Theaetetus* 151d7-e7: It appears ('phainetai') to me that …[so and so]. Because ('gar') of that you must declare ('apophainimenon')…[so and so]. *Politeia* 470a8-b10: Because ('heneka') it appears ('phainetai') to me that…[so and so]…, I declare ('apophainomenou') that… [so and so]. *Politicus* 276e10-a6 [in contrapositive manner]: May we declare ('apophainometha') that …[so and so]? We cannot because it does not yet appear ('oupo phainetai' to me…[so and so]. *Timaeus* 39a4-b2: It appeared ('ephaineto') that …[so and so]; because ('gar')...to appear ('apephainen')... . Similarly, *Protagoras* 340b2-c8 ('phainetai', 'proapophenai,' apephenato'), *Philebus* 67a10-16 ('phanentos', 'pephantai', 'apephenato').
[53] *Theaetetus* translation, Loeb edition, p. 25:'showing that squares …are *[[each, distributively]]* not commensurable in length with the unit of the foot' vs. 'since the number of roots appeared to be *[[collectively]]* infinite'.
[54] [Burn2], p. 266.:'showing that the power of 3 square feet and the power of 5 square feet are *[[each, distributively]]* not commensurable in length with the power of 1 square foot' vs. 'since the powers were turning out to be *[[collectively]]* unlimited in numbers'.Burnyeat is following Levett's translation, even though, as we saw in section 1.4, considers [b] as being thoughts prompted from [a].
[55] [Kn], p. 62-63:'demonstrating that these are *[[each, distributively]]* not commensurable in length with the one foot-power' vs. 'since we recognized the powers to be *[[collectively]]* unlimited in number'.
[56] AST is freely exchanging 'apophainon' with 'ediknuen' 25,42; 'epedeiknuen' 28,42; 'ediknuen' 34,2 and 34,20.

contrary to the Cherniss-Burnyeat thesis, as we shall see.

*(b) The distributive sense of the plural 'hai dunameis' in the crucial statement [b]*

But now the two *plurals*, 'ou summetroi' in [a] and 'apeiroi to plethos' in [b], must be examined under the light thrown by the existence of the bond between [a] and [b]. We note that it makes no sense for the plural 'ou summetroi' used in [a] to have *a* collective meaning (such as that somehow they are all together as a totality incommensurable), its only possible meaning being *distributive*[57], namely that

*'each of the powers* was shown to be incommensurable to the one foot line''

It is then suggested that the plural 'apeiroi to plethos' used in [b], bonding in structure and meaning with [a], does not have a collective (such as that somehow they are all together infinite in multitude), but its meaning must be *distributive* as well:

[b***] 'since *each of the powers* was shown to be infinite in multitude'.

This is then the rendering we propose for [b], on the basis of the foregoing analysis. The same rendering, statement [b***], has been arrived at in section 1.4 by an independent route.

Our conclusion is strengthened by the persistence in using distributive plural for 'hai dunameis' in the sequel 148b1-2:

*'those lines were called powers ('dunameis') that were not commensurable ('ou summetrous') in length to the 'lengths', but commensurable in the squares they can (to the squares that the 'lengths' can).'*

*3.2. Connection between proof of incommensurability by Theodorus and the infinite in the* Anonymous Commentator to Plato's Theaetetus *(AST),36,36-37,29.*

The comments in AST for the crucial statement, included in AST 36,36-37,29, provide unqualified support for the thesis that the plural ('hai dunameis apeiroi to plethos') in the crucial *Theaetetus* 147d7-8 statement is *distributive*, and that *every* power admits of the infinite *separately* in a divisional way.

(a) It will be sufficient to confine ourselves in the brief but revealing sentence AST 36,45-48, appearing as an explanation of the crucial *Theaetetus* statement :

*'because the lines admit of the indefinite ('aoriston') either by increasing ('auxoi') or by dividing ('diairoi') them, [and the lines are limited ('horizontai') by numbers, Theaetetus and his companion passed over to numbers'].*

---

[57]We wish to thank the anonymous referee of an earlier version for bringing our attention to this (standard) terminology. 'A statement in English such as "all squared nonzero real numbers are positive" is called a *distributive plural*. This means that the statement "the square of *x* is positive" is true for every nonzero real number. It can be translated directly into: 'for all x (if x is a non-zero real number, then $x^2$ is positive)'. Not all statements involving plurals in English are distributive plurals. The statement "The agents are surrounding the building" does not imply that Agent James is surrounding the building. This type of statement is called a *collective plural*. Such a statement cannot be translated directly into a statement involving a universal quantifier.' (after 'abstractmath.org').
A statement with a collective plural can be translated directly into a set theoretic statement.

The meaning of the AST statement is clear enough by itself:
the infinity that lines (and powers, considered as lines[58]) admit can be of *only* two kinds, either (a) infinity by division, or (b) infinity by increase,
but it becomes clearer if we compare it with numerous similar statements from ancient authors, both mathematically and philosophically inclined:

Aristotle, *Physics* 233 a24-26
*'For, length ('mekos'), and time, and generally every continuum ('suneches') is said to be infinite ('apeiron') in two senses, either by division ('kata diairesin') or by the endpoints ('tois eschatois').'*

Heron, *Definitiones* 119,1,2-3 (=*Scholia in Eucliden* 5, 4,1-2)
*'Magnitude is increasing ('auxanomenon') and divisible ('temnomenon') ad infinitum'*

Themistius, *In Aristotelis physica paraphrasis* 5,2,187,5-7
*'infinity has a double meaning for length and for every continuum in general*
*either by being divided ('diaireisthai') ad infinitum*
*or by the magnitude having no last point and no limit ('meden eschaton...mede peras')'*

Proclus, *In Eucliden* 184, 17-29
*'...that every magnitude is divisible ('diaireton') ad infinitum is a geometric postulate (184, 17-18)...*
*and that the quantity increases ('auxein') ad infinitum is common to both; for both number and magnitude are subject to it (184, 27-29)'*

Proclus, *in Eucliden* 198,14-15
*'for every continuum is divisible ('diaireton') and increasable ('auxeton') ad infinitum'*

Simplicius, *Commentary to Aristotle's Physics* 467,8-12
*'the mathematical magnitudes are divided ('diaireisthai') and increased ('auxesthai') ad infinitum'*

The comparison of the AST statement 36,45-48 with these statements, and the obvious meaning they all have, makes it clear that the plural lines ('ai grammai') in the AST statement in 36, 45-48 is *distributive*, since both kinds of infinity are about a *single initial given line*; the line admits of the infinite either because it is divided into an infinite multitude of parts, or because it is increased by adding to it an infinite multitude of parts. That a given line is divisible (into two lines, hence, recursively, *ad infinitum*) follows essentially from Propositions I.3, 10 of the *Elements*, and that a given line is increasable *ad infinitum* is the content of the second postulate of the *Elements*.

(b) But now we must also bear in mind that the AST statement in 36, 45-48 is meant to be in explanation of the crucial statement *Theaetetus* 147d7-8, reproduced just above in lines 36, 37-40 of AST.
Thus the crucial statement:
[b] 'because the powers were shown to be infinite in multitude'
 is read, after the explanation supplied by AST 36,45-48, in a more precise, if forced, way as:
[b1] 'because the powers were shown to be infinite in multitude either by division or by increase'.

---

[58] Cf. footnote 34

But now the plural appearing in the statement ceases to be ambiguous and cannot but be a *distributive* plural and thus mean:

[b2] 'because *each power* was shown to be infinite in multitude, either by division or by increase'.

This is the third, and most explicit, indication that the plural in the crucial statement is a distributive one.

But there is more. A power a
(i) is according to its definition *a fixed finite line* (satisfying $a^2=Nr^2$ for some non-square number N) that does not in any way increase, and
(ii) is infinite in multitude, according to the *Theaetetus* 147d7-8 crucial statement, just analysed.
Therefore a power cannot possibly admit of the infinite by increasing *ad infinitum*, because then the line would not be fixed and finite (as in (i)). Thus the only possibility left is for the crucial statement is to mean that

[b3] 'because *each power* was shown to be infinite in multitude by division'.

Finally the question arises: how was such a thing *shown*? Here the link between 'ephainonto' and the 'apephenen', suggested earlier, provides the only reasonable explanation:

[a]+[b3] Theodorus showed ('apephenen') that each power is not commensurable in length to the one foot line, because each power was shown ('ephainonto') to be infinite in multitude by division.

In turn this can only imply that Theodorus method was anthyphairetic, that Theodorus was employing in some form *Proposition X.2* of Euclid's *Elements*, according to which if two unequal magnitudes have an infinite anthyphairesis, then these two magnitudes are incommensurable.

*3.3. AST 37,23-26: 'epei apeiroi ephainonto hai kata meke dunameis'*

There are two more versions of the crucial statement in AST 25,40-26,13 and 37, 23-29. We will consider here the AST 37,23-29 passage, and argue that and it confirms both the distributivity of the plural and the connection of infinity and incommensurability of each power separately. We will deal with the AST 25,40-26,13 in section 4.2.

[AST 37,23-29] *'because the powers according to length ('hai kata meke dunameis') were shown to be infinite ('apeiroi') [[Theaetetus andhis companion]] attempted to circumscribe ('perilabein') all these powers by means of a common name.'*

This version ('epei apeiroi ephainonto *hai kata meke dunameis'*, 37,23-26) differs from the original at two points: (i) the original 'apeiroi to plethos' is shortened to 'apeiroi', a change of convenience; and (ii) the original 'hai dunameis' is expanded to 'hai kata meke dunameis'. The meaning of this addition (ii) seems to be the following: powers can be considered in two ways, as lengths-sides, linearly, and as squares[59]; infinity of powers was shown by Theodorus *only* as lengths ('kata meke'), and not as squares. It follows that expansion (ii) makes no sense in the Burnyeat interpretation, as all powers taken together, collectively, form an infinite set, irrespectively of whether they are considered as lines or as squares. The appearance of infinity for powers as lines and not for powers as squares makes good sense only if we interpret the plural as distributive and the infinity as divisional-anthyphairetic.

*4. The crucial statement: connection with Theaetetus and his companion's method*

*4.1. The connection between the Theaetetus and his companion's idea and method and the infinite in the* Theaetetus *147d-e passage.*

In section 3 we have examined the meaning of the crucial statement [b] by exploiting its connection with Theodorus' lesson and demonstrations. But the crucial statement [b] is transitive from Theodorus' lesson to Theaetetus and his companion's idea and method, and the meaning of the crucial statement [b] can also be examined by exploiting its connection with Theaetetus and his companion's idea and method.

In the *Theaetetus* 147d8-148b2 passage the transitive crucial statement [b], 'since the powers were shown to be infinite in multitude', is followed by the statement that Theaetetus and his companion had the idea to collect in one ('sullabein eis hen'), later described as 'circumscribe' ('perilambanein', 148d6, the powers (147d8). For this goal, they divided the set of all numbers in square and non-square numbers (147e5-148a4), and they defined line a, such that $a^2=Cr^2$, where r is the one foot line, *length* if C is a square number and *power* if C is non-square. Thus in a power there is incommensurability of sides and commensurability of squares (148a6-b2)[60]. This approach is judged as quite successful by Socrates, to the point that it must be taken as model of imitation for his philosophical quest on the meaning of knowledge (148c-d). However from Plato's account, just outlined, the mathematical method and contribution by Theaetetus and his companion is by no means clear[61], although it is generally related to Proposition X.9 of the *Elements*.

As we will see, the Anonymous Commentator AST throws substantial light on their method. Our interest in the present paper is limited to the light thrown in relation to the crucial statement [b], and we will be content with some heuristic comments for anything beyond that goal.

(a) In general *infinity is circumscribed by means of a whole entity*

According to AST 37,3-12, it is a general mathematical and philosophic principle that whenever an infinity occurs, there must be an attempt to collect ('sullabein' *Theaetetus* 147d8) - circumscribe ('perilambanein', *Theaetetus* 148d6; AST 26,11; 37,5; 37,10; 37,29; 37,43; 37,46; 45,48; 46,39) - define ('horizein', AST 37,1;37,11) - limit ('periorizein', AST 42,32) - subjugate ('hupotaxai', AST 26,12) the infinity in question.

---

[59] Cf. footnote 34
[60] Cf. footnote 34
[61] {Burn1], p. 505-509

[AST 37,3-12] *'because the infinite ('apeiron') is uncircumscribed ('aperilepton') and indefinite ('aoriston'), intelligible thought ('dianoia') must, to the extent possible, circumscribe ('perilambanein') and define ('horizein') it by means of a whole entity ('katholiko tini').*

Thus in order to circumscribe an infinite a one and whole entity must be found, if possible.

(b) *Theaetetus and his companion achieved the circumscription of the infinity in powers (147d7-8) by passing over to numbers (147d8-148b2)*

In order to achieve the desired circumscription of the specific infinite of powers that arises in the *Theaetetus* 147d7-8 crucial passage, Theaetetus and his companion passed over to numbers:

[AST 36,45-37,3] *'because the lines admit of the indefinite ('aoriston') either by increase ('auxoi') or by division ('diairoi'), and the lines are limited ('horizontai') by numbers, Theaetetus and his companion passed over to numbers'.*

(c) *The three steps in Theaetetus and his companion' s method*

Thus we are informed that Theaetetus and his companion in their attempt to circumscribe the infinite in powers (i) passed over (from the infinite) to numbers (b), and (iii) the infinite was indeed circumscribed by means of an entity described as One and Whole. There is thus an imlied missing link (ii) from the numbers to the One and Whole entity.
We thus concude that Theaetetus and his companion achieved the circumscription of the infinity of a power in three steps:
(i) from the infinite of a power they passed over to numbers (by (b));
(ii) the passage over to numbers allows for the formation of an entity which is described as one/whole ('katholiko tini')/logos and whose name serves as the name of all (the infinite parts); and,
(iii) by means of the one and whole entity the infinite is circumscribed (by (a)).

Step (i) is explained in AST in unambiguous terms, as we describe in (d). For the purpose of the present work it is not necessary to commit ourselves as to (and subsequently attempt to establish) the precise meaning of the entity described as one/whole/logos (step (ii)) and of the circumscription and definition of the specific infinite that results in the *Theaetetus* 147d7-148b2 passage (step (iii)) and is meant to give an account of the important mathematical achievement of Theaetetus and his companion. But some informal remarks, which will not be used for drawing conclusions about our sole interest in this paper, which is the use of infinity in the crucial statement [b], may be helpful and hopefully illuminating and will be given in (e).

(d) *The nature of step (i)*

[1] By [AST 44,50-45,5] the passage over to numbers [step (i)] is more carefully expressed as the passage from magnitudes to numbers.

[AST 44,50-45,5] *'Theaetetus and young Socrates passed over from the less clear to the more clear, and to the more whole ('ta katholikotera') because they passed over **from magnitudes** to numbers.'*

[2] By [AST 26,13-18] by the expression 'the passage over to numbers' [in step (i)] the expression 'the passage to a ratio 'as number to number'', a ratio that is always commensurable' is meant.

[AST 26,13-18] *'thus they came ('elthon')[62] to number as a consequence ('dia to akolouthon ?'[63]) of the fact that all the numbers are commensurable to each other'*,

We may infer from [1] and [2] that
[3] the passage over to numbers [step (i)] is really the passage from a ratio of magnitude to magnitude to a ratio described 'as number to number', equivalently as a commensurable ratio.

Expanding, in AST 41,17-45,5 on the quite brief *Theaetetus* 148b2 passage, to a power-length approach for cubes similar to that of squares the Anonymous Commentator explains:

[AST 42,13-33] *'for the cube itself to the cube has commensurable ratio, as number to number, while the sides are incommensurable ,,, thus ('oun') for the solids as well they passed over to numbers in order to set a limit ('periorisosi') by means of a whole entity'*,

We infer from [AST 42,13-33] and [3], (modifying it an obvious way from cubes to squares and) making use of the propositions, stated and proved as Propositions 5, 6 of Book X of the *Elements* stating that the ratio of two magnitudes is commensurable if and only if that ratio is as number to number, that

[4] for the power defined by the equality $a^2=Cr^2$, where r is the one foot line and C a non-square number, the commensurable as number to number ratio is the commensurable ratio $a^2$ to $r^2$, equal to the number to number ratio C to 1, and the magnitude to magnitude ratio is the incommensurable ratio a to r.

We infer from [3] and [4] that

[5] step (i), namely the passage over to numbers is the passage from the incommensurable [magnitude to magnitude] ratio a to r to the commensurable as number to number ratio $a^2$ to $r^2$, namely the passage obtained by **squaring** from the incommensurable ratio of the sides a to r to the commensurable ratio as the number C to number 1 ratio of the squares $a^2$ to $r^2$.

At this point we take into account (b) above, namely that
[6] by [AST 36,45-37,3 & 37,3-12] the passage over to numbers [step (i)] was undertaken by Theaetetus and his companion for the purpose of circumscribing and limiting the infinite in powers [step (iii)].

Combining [5] and [6] we see that

---

[62] Cf. AST 42,3 for accepting 'elthon' as equivalent to 'metebesan'
[63] Cf. AST 47,46-47 and 63,48 for similar expressions.

[7] step (i), on the one hand, by [5], is a step from incommensurability to commensurability, while on the other hand, by [6], is a step from the infinity of powers to a condition that prepares the ground for limiting and circumscribing the infinite.

Thus the infinity of powers must be closely related to incommensurability of powers to the one foot line, while commensurability of their squares must be instrumental for the circumscription and limiting of the infinite. But the only conceivable relation of commensurability and as number to number ratio to the finite is the known proposition that *finite anthyphairesis* characterises commensurability; and the only possible relation of incommensurability to the infinite is the known proposition that *anthyphairetic division ad infinitum* characterises incommensurability (Propositions X.2 & 3 of the *Elements* and the discussion in section 1.3). We conclude that

[8] the infinity of the powers in the crucial sentence [b] in the *Theaetetus* 147d7-8 is the infinity that passed over to numbers, and is thus the infinity of anthyphairesis resulting from the incommensurability of the power with respect to the one foot line, and

[9] the passage to numbers is the passage by squaring to the commensurability, and resulting finite anthyphairesis, of the square $a^2$ with respect to the one foot square $r^2$, a passage from the infinite to the finite [step (i)] meant to produce an entity described as One and Whole [step (ii)], by means of which the initial antyphairetic infinite of the power is circumscribed and limited [step (iii)].

It is clear from [8] (and [9]), clarifying step (i) of Theaetetus and his companion contribution, that the AST commentary on the Theaetetus and his companion idea and mathematical discovery that the crucial statement [b] in *Theaetetus* 147d7-8 is in opposition to the Burnyeat interpretation, supporting the distributive plural, and is in opposition to Cherniss Byurnyeat neutrality thesis, supporting the anthyphairetic interpretation.

(e) *Heuristics on the nature of steps (ii) and (iii)*

The process described as circumscribing the infinite does not mean that the original infinite ceases to be infinite, on the contrary it remains infinite, a circumscribed and limited infinite.[64] Circumscription is a circular containment with a strong connotation of circularity and periodicity[65]. So in the case of the infinity-incommensurability of a power, which as mentioned in (a) is anthyphairetic, a reasonable interpretation of the circumscription of infinity (step (iii)) would be a circular, periodic anthyphairesis. This is in full agreement with the mathematics as we know it, since the anthyphairesis of every power is indeed known to be periodic[66]. Periodicity of anthyphairesis is closely related to a deep and seemingly paradoxical sense of Oneness, namely that of the Whole that cannot be destroyed but remains whole after any division, namely an entity in which every part of the whole is the same as the whole, a description that fits with the 'whole' ('holon') and 'one', as described in the *Parmenides* 142d-e, 144b-145a, *Sophistes* 245a, and *Philebus* 15a-b. Thus step (iii) may reasonably be said to be caused by an entity that may be described as One and Whole (step (ii)).

---

[64] As is made clear e.g. in Plotinus, *Enneades* 2,4,16; 6,6,3,9-16, in the infinite of Division and Collection, and Iamblichus, *De Communi Mathematica Scientia* 3.

[65]Cf. a non-exhaustive sample: *Republic* 546b4 'periodos…perilambanei' (cf. Proclus, *eis Politeian* 2,30,10-14 adds 'anakukleisthai…apokathistamenon'); *Sophistes* 235a-c , 'perilambanein' assumes an expressly circular character of containment: 'we have got him into a kind of encircling net ('auton perieilephamen en amphiblestriko tini', 235b1); *Timaeus* 32c5-34a7 the sphere is 'perieilephos' all the other figures. The same is true in Euclid, *Elements* Definitions XI. 14,18, 21and Propositions XIII.13-17. Plotinus, *Enneades* 2,2,1,11 'kuklo'…'perilambanein'; 5,1,8,18-22 'sphairas'…'perieilemmena'; 6,4,2,35-36 'perilabein' associated with 'perithein'; 6,4,2,35-36'kuklo'…'perilabein';6,8,18 (similar); 6,8,11,22-37 (similar). Iamblichus, *De Mysteriis* 1,9,31 'kuklo perilabein'. Porphyrius, *In Aristotelis categorias* 4,1,120,16 'kuklo perilabein'. Galenus, *In Hippocratis librum de officinamedici commentarii* 18b, 712,14 'kuklo…perilabein'.Proclus, *eis Timaion* 1,209,26-'kuklo perieileptai';

Stated in mathematical language, [9] tells that Theaetetus & his companion achieved the circumscription of the infinite for powers, namely the periodicity of the anthyphairesis of each power, infinite because of incommensurability, because in powers it is possible to pass over by squaring from the incommensurability-infinity of the power to the commensurability-finiteness. Thus, according to the *Theaetetus* account in 147d7-148b28 and the AST commentary on that account, Theaetetus and his companion discovered and proved the following fundamental

*Proposition.* If a is a power (namely, if $a^2=Cr^2$ for some non-square number C, where r is the one foot line),
then
(i) a is incommensurable to r, and
(ii) the anthyphairesis of a to r is periodic by means of the logos criterion.

The statement of this proposition is closely related to Proposition X.9 of the *Elements*[67].

*4.2. AST 26,6-8: 'epei toinun apeira en ta toiauta tetragona'*

We finally examine the AST 25,40-26,13 passage, containing a version of the crucial statement [b]. In it there is no mention of the term 'power', but only of 'square'[68]; it starts with the phrase corresponding to the *Theaetetus* 147d3-6,

*'Theodorus was proving ('edeiknuen') to those around Theaetetus that the three-foot and the five-foot, is incommensurable [to the one foot square] according to the sides ('kata tas pleuras'), from which each was produced. And numbering ('exarithmoumenos') the incommensurable squares '(ta asummetra tetragona') came till as far as the 17-foot square.'*

Now AST comes to a version of the crucial statement [b], involving squares:

*'Because these ('toiauta') squares were infinite ('apeira'), those around Theaetetus attempted to circumscribe by means of some whole entity, so as to subjugate by means of a single name'.*

---

1,247,16-248,6 'perilambanein…peritheon…peri kentron…perichoreuontos autou kuklo'; 2,71,16-72,5 'perileptikon'…'sphairikon'; 2,109,1-; 3,272,22- 'kata tina periodon'…'he periodike poiesis'…'perilabein'; *In Parmeniden* 807,24 'perilabein' … 'perithei' … 'perichoreuei kuklo'. Syrianos, *eis Hermogenen*, 11,2 'kuklo perilambanein'. Simplicius, *In Aristotelis physicorum libros commentaria* 9,62,19-20 'perilabon kuklo'.

[66] In other words, the continued fraction of a quadratic irrational is (eventually) periodic (cf. [F], 9.1).

[67] As pointed out in *Scholion in Eucliden* X.62, the hypothesis of Proposition X.9 is more general than that of the Proposition implied in the *Theaetetus* text. Cf. [Burn1], p. 505-509.

[68] Cf. footnote 34

The 'toiauta' squares in 26,7 refers to the incommensurable ('asummetra') squares, as in 26,2-3, with the understanding, as in 25,48-49, that they are incommensurable according to the sides ('kata tas pleuras'), and thus the phrase
'epei toinun apeira en ta *toiauta* tetragona' (26,6-8)
should be read as
'epei toinun *apeira* en ta *asummetra kata tas pleuras* tetragona'.

The statement, considered by itself, is ambiguous: it is not clear if 'ta tetragona' is a distributive or collective plural. The similarity of this statement with the statement

*'because powers according to lengths were infinite'* (AST 37, 23-26),

examined in section 3.3, in which the corresponding plural was found to be distributive, might suggest that the plural here is distributive too. But the stronger argument for distributivity is the fact that the statement in question [AST 25,40-26,13] is followed by the statement [AST 26,13-18], which, in section 4.1 (d) [2] above, was shown to be a passage from the infinite over to numbers, namely to commensurability and finiteness, and thus a passage from the incommensurable infinite, an infinite that makes sense only for each power or square separately. This shows the distributive and anthyphairetic interpretation of the statement in question.

*5. Conclusion*

Scholars have suggested various reconstructions for the incommensurability proofs of Theodorus' lesson in the *Theaetetus* 147d3-e1 passage, either anthyphairetic (notably by Zeuthen, Becker, van der Waerden) or non-anthyphairetic (notably by Hardy & Wright and Knorr), but the dominant view has been the thesis by Cherniss and Burnyeat on the neutrality of the platonic text in relation to any of these reconstructions.

The weak points of Knorr's reconstruction, according to which Theodorus ended his demonstrations before the case C=17, were discussed in section 1.1, while the argument by Knorr on the impossibility of a rigorous proof of Proposition X.2 of the *Elements* was countered in section 1.3, with the result that the possibility of an anthyphairetic method for Theodorus' demonstrations not to be excluded.

The difficulty with Burnyeat's interpretation of the crucial sentence in the *Theaetetus* 147d7-8 was discussed in section 1.4, resulting in a problem for the Cherniss-Burnyeat neutrality thesis.

The anthyphairetic nature of the correlation of incommensurability with either infinite divisibility or a principle of the infinite is discussed in section 2.1. The anthyphairetic interpretation of Theodorus method of proof is based on the argument that the crucial statement 'because the powers were shown to be infinite in multitude' in 147d7-8 is interpreted according to the existing evidence, only with the adoption of a distributive plural, contrary to the Burnyeat interpretation, resulting for an anthyphairetic method for Theodorus' proofs of incommensurability, contrary to the Cherniss-Burnyeat neutrality thesis. The distributive and anthyphairetic approach is supported not only by the *Theaetetus* 147d3-e1 passage itself (in 3.1), but also by the relevant comments in *Scholia in Platonem* (2.2) and extensively in the *Anonymi Commentarius In Platonis Theaetetum* (3.2, 3.3 and 4).

On the basis of these arguments we conclude (a) that the only credible non-anthyphairetic reconstruction, namely the one by Hardy and Wright, is excluded, because it is inconsistent with the distributive interpretation of the crucial statement, and (b) that contrary to the neutrality thesis of Cherniss and Burnyeat, Plato does inform the reader on the nature of the method of Theodorus proof of incommensurabilities and that Theodorus' method of proof is anthyphairetic.

*Department of Mathematics,*
*Athens University, Athens 157 84, Greece*
*snegrep@math.uoa.gr*

*Bibliography*

[AST] *Anonymi Commentarius In Platonis Theaetetum.*
[SP] *Scholia in Platonem*

[Be] O. Becker, *Eudoxos-Studien I. Eine voreudoxische Proportionlehre und ihre Spuren bei Aristotles und Euklid*, Quellen und Studien zur Geschichte der Mathematik, Astronomie und Physik, Abteilung: Studien 2 (1933) 311-333.
[Burn1] M. F. Burnyeat, *The Philosophical Sense in Theaetetus' Mathematics*, Isis 69 (1978), 489-513.
[Burn2] M. F. Burnyeat, *The Theaetetus of Plato; with a translation of Plato's Theaetetus by M.J. Levett, revised by M. Burnyeat*, Hackett Publishing Company, Indianapolis, 1990.
[C] H. Cherniss, *Plato as Mathematician*, The Review of Metaphysics 4 (1951) 395-425.
[F] D. Fowler, *The Mathematics of Plato's Academy, a new reconstruction*, Second edition, Clarendon Press, Oxford, 1999.
[vF] K. von Fritz, *Theodorus*, in Pauly Wissowa *Realencyclopaedie der Classischen Altertumwissenschaft* 10 (2), cols. 1811-1825, 1934.
[Ha] R. Hackforth, *Notes on Plato's Theaetetus*, Mnemosyne 10 (4) (1957) 128-140.
[HW] G. H. Hardy and E.M.Wright, *Introduction to the Theory of Numbers*, Oxford, Clarendon Press, 1938.
[He] T. L. Heath, *The Thirteen Books of Euclid's Elements*, Second edition, Three volumes, Cambridge University Press, Cambridge, 1926.
[Ka] J.-P. Kahane, *La Theorie de Theodore des corps quadratiques reels*, L' Enseignement Mathematique 31 (1985) 85-92.
[Kn] W. R. Knorr, *The Evolution of Euclidean Elements: A Study of the Theory of Incommensurable Magnitudes and Its Significance for Early Greek Geometry*, Reidel, Dordrecht, 1975.
[Mo] G.R. Morrow, *Proclus: A Commentary on the First Book of Euclid's Elements, Translated, Introduction, Notes*, Princeton University Press, Princeton, 1970.
[N1] S. Negrepontis, *The Anthyphairetic Nature of the Platonic Principles of Infinite and Finite*, in Proceedings of the 4th Mediterranean Conference on Mathematics Education, January 2005, Palermo, Italy, 3-26.
[N2] S. Negrepontis, *The Anthyphairetic Nature* of *the Geometry and the Philosophy of the Pythagoreans*, in volume edited by D. Anapolitanos, *Instances and Durations,, Thirteen Philosophical Essays*, Nephele, Athens, 2009, p. 191-282 [in Greek].
[N3] S. Negrepontis, *Plato's Theory of Knowledge in the* Sophistes, Proceedings of the Colloquium on 'La Demonstration de l' Antiquite a l'age classique', June 3-6, 2008, Paris, to appear.
[Ta1] A. E. Taylor, *Forms and Numbers: A study in Platonic Metaphysics*, Mind (I) 35 (1926) 419-440; (II) 36 (1927) 12-33.
[Ta2] A. E. Taylor, *Plato, The Man and his Work*, Methuen, London, 1926.
[Thoms] W. Thomson, *The Commentary of Pappus on Book X. of Euclid's Elements. Arabic Text and Translation. With Introductory Remarks, Notes, and a Glossary of Technical Terms by Gustav Junge and William Thomson*, Harvard University Press, Cambridge, 1930.
[Thomp] D'A. W. Thompson, *III. Excess and Defect: or the Little More and the Little Less,* Mind 38 (1929) 43-55
[vdW] B. L. van der Waerden, *Science Awakening*, translation by A. Dresden of *Ontwakende Wetenschap* (1950), Noordhoff, Groningen, 1954.

[Ze] H. G. Zeuthen, *Sur la constitution des livres arithmetiques des Elements d' Euclide et leur rapport a la question de l' irrationalite*, Oversigt over det Kgl. Danske Videnskabernes Selskabs Forhandlinger 5 (1910) 395-435.